\documentclass[10pt]{amsart}
\usepackage{amssymb, amsmath, latexsym}
\usepackage{graphicx}
\usepackage[dvipsnames,svgnames,x11names,hyperref]{xcolor}
\usepackage[a-1b]{pdfx}
\hypersetup{colorlinks,breaklinks,
             urlcolor=blue,
             linkcolor=blue}

\newtheorem{theorem}{Theorem}[section]

\theoremstyle{definition}

\theoremstyle{remark}

\numberwithin{equation}{section}
\def\F{\mathbb{F}}

\def\Z{\mathbb{Z}}
\def\Fq{{\mathbb F}_q}

\def\PP{{\mathbb P}}

\def\hp3{\widehat{\mathbb P}^3}
\newcommand{\Ev}{\mathrm{Ev}}
\newcommand{\RM}{\mathrm{RM}}
\newcommand{\PRM}{\mathrm{PRM}}

\def\Aff{{\mathbb A}}
\begin{document}

\title[AGCT: Homage to Gilles Lachaud]{Arithmetic, Geometry, and Coding Theory:  Homage to Gilles Lachaud}

\author[Ghorpade]{Sudhir R. Ghorpade}
\thanks{Sudhir Ghorpade is partially supported by DST-RCN grant INT/NOR/RCN/ICT/P-03/2018 from the Dept. of Science \& Technology, Govt. of India, MATRICS grant MTR/2018/000369 from the Science \& Engg. Research Board, and IRCC award grant 12IRAWD009 from IIT Bombay.}
\address{Department of Mathematics, 
Indian Institute of Technology Bombay,\newline \indent
Powai, Mumbai 400076, India}
\email{\href{mailto:srg@math.iitb.ac.in}{srg@math.iitb.ac.in}}

\author[Ritzenthaler]{Christophe Ritzenthaler}
\address{Univ Rennes, IRMAR, 
  Campus de Beaulieu, %
  35042 Rennes, %
  France}
\email{\href{mailto:christophe.ritzenthaler@univ-rennes1.fr}{christophe.ritzenthaler@univ-rennes1.fr}}

\author[Rodier]{Fran\c cois Rodier}
\address{Aix Marseille Universit\'e, CNRS, Centrale Marseille, \newline \indent
Institut de Math\'ematiques de Marseille, UMR 7373, 13288 Marseille, France}
\email{\href{mailto:francois.rodier@univ-amu.fr}{francois.rodier@univ-amu.fr}}

\author[Tsfasman]{Michael A. Tsfasman}
\thanks{Michael Tsfasman is supported in part by ANR project FLAIR (ANR-17-CE40-0012).}
 
\address{CNRS, Laboratoire de Mathematiques de Versailles (UMR 8100), France\newline \indent
Institute for Information Transmission Problems, Moscow, Russia \newline \indent
Independent University of Moscow, Russia}
\email{\href{mailto:mtsfasman@yandex.ru}{mtsfasman@yandex.ru}}
\date{}
\begin{abstract}
We give an overview of several of the mathematical works of Gilles Lachaud and provide a historical context. This is interspersed with some personal anecdotes highlighting many facets of his personality. 
\end{abstract}
\maketitle

The AG${\rm C}^2$T Conference at CIRM\footnote{CIRM = Centre International de Rencontres Math\'ematiques = International Centre for Mathematics Meetings (\href{https://www.cirm-math.fr}{https://www.cirm-math.fr}) is based in Luminy, near Marseille, France.} 
in June 2019 was dedicated to 
Gilles Lachaud (26 July 1946 -- 21 February 2018) who was a founder of this series of conferences which began in 1987. Two speakers spoke exclusively about the work and persona of Gilles Lachaud, while many others paid rich tributes to him during their talks\footnote{The slides of most of the talks in AG${\rm C}^2$T-2019, and videos of some of them, are available at: \href{https://conferences.cirm-math.fr/1921.html}{https://conferences.cirm-math.fr/1921.html}}. This article is a homage to Gilles and it attempts to give an outline of his mathematical journey and glimpses of his work and personality. It is divided in four parts. The first part, written by Rodier, gives an account of Gilles' early works including his Ph.D. thesis related to automorphic forms, works on Warings problem and real quadratic fields, and some of his earlier articles on coding theory. Next, Tsfasman recounts the beginnings of algebraic geometric codes,  his interactions with Gilles Lachaud, and some of their joint work. The third part, by Ghorpade, describes some of Gilles' work on  linear codes and on counting the number of rational points of varieties over finite fields, as well as some related developments. Finally, Ritzenthaler gives a personalized account of Gilles' work that deals mainly with non-hyperelliptic curves and abelian varieties. 

\smallskip

\section{Gilles Lachaud's early works} 
\medskip
\centerline{Fran\c cois Rodier}
\medskip
I met  Gilles Lachaud in 1969 when he had just been appointed to what was still called the University of Paris, as Roger Godement's assistant. He was a pure product of the university. He did not want to go to the so-called ``Grandes \'Ecoles". 
Gilles had always refused this royal road.

During the  academic year 1969-70, Roger Godement had {left} Paris for a six-months stay at the Institute for Advanced Study in Princeton, USA, to work with, among others, Robert Langlands. Meanwhile, Godement's assistants organized the same DEA (Masters degree) course that he offered each year. With Gilles, there were Miss Chambon, Jean-Pierre Labesse, and Michel Duflo, assisted by Fran\c cois Bruhat. It was by attending these classes on {\sl Representations of Locally Compact Groups} that I first met Gilles. He had lectured on spectral theory, a field in which he specialized.

Having started a thesis myself under the direction of Godement, I met Gilles quite often, especially during the weekly sessions of seminars. He  came to the symposium in 1971 in Budapest where we had the opportunity to cross the Iron Curtain to meet the Soviet School of  Group Theory, and first Israel Gelfand, who had been authorized to cross the border of the USSR to come to Hungary. Then, in the evening, Gilles organized some parties with the younger people in the conference.

In 1972, an International Summer School was organized by the University of Antwerp, on  {\sl Modular functions of one variable}. But it was financed by NATO. The dislike of the army by Godement which was  well known, made him prefer a symposium organized by the AMS in Williamston, in the United States on {\sl Harmonic analysis on homogeneous spaces}. Gilles was one of the mathematicians who followed Godement instead of attending the symposium in Antwerp. This was the first time that I did my presentation in English in front of an amphitheater filled with people that I considered more capable than me. So Gilles offered me to make me rehearse it.
Gilles was more sure of himself in English, but he still asked me to assist in exchange for his own rehearsal.

After publishing two articles in French in the journal Inventiones:
{\sl Analyse spectrale des formes automorphes et s\'eries d'Eisenstein} \cite{rod-lac1} in 1978, and {\sl Variations sur un th\`eme de Mahler} \cite{rod-lac2} in 1979, Gilles defended his  thesis of  ``Doctorat d'\'Etat" at Paris Diderot University - Paris 7 in June 1979 on  {\sl Analyse spectrale et prolongement analytique: S\'eries d'Eisenstein, Fonctions Z\^eta et nombre de solutions d'\'equations diophantiennes} \cite{rod-lac6} under the direction of Roger Godement. 
In this thesis, Gilles reproves the results of Selberg 
on the meromorphic prolongation of the Eisenstein series by a method of perturbing the continuous spectra of self-adjoint operators. He also obtains the analytical extension of a type of generalized zeta functions, and applies these results to diophantine equations.
The Rivoire prize from the University of Clermont-Ferrand was awarded to him in 1979 for his thesis.

More specifically, in the article {\sl Spectral analysis of automorphic forms and Eisenstein series}, 
Gilles is concerned with the spectral decomposition of the left regular
representation $\pi$ of $G$ on $K\backslash G/\Gamma$ for $G$ a real, linear, connected simple Lie group of
rank 1, $K$ a maximal compact subgroup, $\Gamma$ a discrete subgroup such that $G/\Gamma$  has
finite invariant volume. 
Then, by a general principle of functional analysis, there exists a spectral measure $m$ on the spectrum of the algebra $L^1 (G, K)$ of 
 complex integrable functions on $G$  bi-invariant under $K$. This spectrum  is a closed subset of the complex plane, such that  we have
$$
(\pi (F) f, g) = \int \tilde F (\lambda) dm_{f,g}(\lambda)
$$
by noting $\tilde F$ the Gel'fand transform of the function $F\in L^1 (G, K)$, and $dm_{f,g}$ the complex measure associated with the 
spectral measure $m$ and the elements $f, g$ of $L^2 (K \backslash G /\Gamma)$.
Some eigenfunctions of $\pi (F)$ are given by Eisenstein's series. In the case of
$G = SL(2,\mathbb{R})$, $\Gamma = SL(2,\mathbb{Z})$,  $K= SO(2)$ they are written
$$E_s(z) = {1\over2} y^s \sum_{(c,d)=1} |cz+d|^{-2s}, \qquad \Re e\ s > 1, z = x+iy, y > 0,$$
since one can identify the symmetric space $K\backslash G$ with the upper half-plane. 
 They extend analytically to the entire complex plane, and their value on the spectrum of the algebra $L^1 (G, K)$ allows to calculate the continuous part  of the 
 measures $dm_{f, g} $.
It was A. Selberg who stated for the first time these results   in the case the group $G$ is of  rank one on the field of rationals, following his own method.

The aim of Gilles is to use another method based on the theory of perturbations of the continuous spectra of self-adjoint operators. It consist in finding a representation $\pi'$ 
of $G$ in  $L^2 (K \backslash G /\Gamma)$ fulfilling two properties. First its spectral analysis is   known, second
 it is
close to the representation $\pi$ in a suitable sense, namely that the difference $\pi (F) - \pi '(F)$ must satisfy a certain compactness criterion for  sufficiently numerous $F\in L^1 (G; K)$.
This method, known from Selberg and Gel'fand, was developed by Faddeev in the case of the group $SL (2, R)$. 
Hence, the ideas used by Gilles in proving his main
theorems are not completely different from known ones. But, after using many general
theorems, the theory has been given a smooth and systematic look.
%

In the article {\sl Variations on a theme of Mahler}, 
 he studied the function
$$ N (t) = \# \{x \in \mathbb{Z}^n | |F (x)| \le t \} $$
 for multivariate polynomials $ F $ with $n$ variables on   $\mathbb{Z}$.
More precisely let $S = \{p_1, \dots, p_s\}$ be a finite set of prime numbers,  $P = p_1 \dots p_s$ and $\mathbb{Z}^n (S)$ 
the set of vectors $x =  (x_1,\dots, x_n)$ of $\mathbb{Z}^n$ whose coordinates are prime to $P$ as a whole.
For any prime $p$, we denote $|x|_p$ the absolute value of $x \in \mathbb{Q}$ associated with the place defined by $p$. 
Denote $|x|_0$ the Archimedean absolute value of $\mathbb{Q}$. If $x\in \mathbb{Q}^n$, we write:
$$F_S (x) = | F (x) |_0 \prod_{p\in S} | F (x) |_p.$$

Gilles proves the following theorem about $ N_S (t) = \# \{x \in \mathbb{Z}^n | |F_S (x)| \le t \} $.
 \begin{theorem}
With the previous notations, let us suppose the forme $F$ an\-isotropic on $\mathbb R$ and on $\mathbb  Q_p$ when $p\in S$. Then
$$N_S (t) = V_S t^{n/d} + O (t ^{(n- \theta) / d}),$$
where $V_S$ is a constant (which he computes) and where $\theta> 1$.
\end{theorem}
He deduced the analytic continuation of certain  Zeta functions.

Strangely enough, because the theory of automorphic forms has developed strongly, under the impetus of Langlands among others, these publications have had little echo. But Gilles's interest in this theme has not disappeared since, as we can see in several of his later publications, and perhaps one of the last  {\sl The distribution of the trace in the compact group of type $G_2$} \cite{rit-lac6}, where he used the theory of representations of Lie groups.
This article responds to a request from J.-P. Serre concerning the calculation of a certain exponential sum of degree 7, involving the $ G_2 $ group as
in the theory of Sato-Tate.

Gilles was interested in other things than the theory of automorphic forms. He was also concerned about number theory and history of mathematics.

He  gave  lectures at various universities about some points of history of mathematics, the foundations of arithmetic and the connections between mathematics and philosophy, mysticism and magic. In the article {\sl Nombres et arch\'etypes: du vieux et du neuf [Number and
archetypes: something old, something new]} \cite{rod-lac3}, he suggested: mathematics and philosophy have been estranged since the time of Diophantos, whose number theory can be viewed  as a succession of relations between objects whose
status is determined axiomatically and mathematical truth is determined internally, without reference to external reality;
there is another tradition of mathematics working under mystical influences, e.g., Euclid's Pythagorean emphasis on perfect numbers and solids. This led Gilles to speculate on the connections
between mathematics and magic with the philosophers like Plato, Iamblichus and Ferdinand Gonseth.

He  contributed to the publication of  {\sl Les Arithm\'etiques} by Diophantus in 1984 ({\sl R. Rashed, Diophante. Les Arithm\'etiques, t. Ill : Livre IV, ccvi + 162 p. ;  Paris, Les Belles Lettres, 1984}).
In 
 \cite[p. 367]{rod-cma},  Karine Chemla writes about this book:
``R. Rashed in the publication of third book of Diophantus makes in his commentary, after having gathered in an abstract of algebraic geometry, written in conjunction with the number theorist Lachaud, the set of techniques that the translation of the text of Diophantus supposes. This abstract thus presents the double virtue of being a well-made mathematical introduction and of being adapted to the text of Diophantus. We have here a most remarkable example of what the combination of mathematics and history can produce."
Indeed, Gilles wrote about fifty pages of introduction to modern algebraic geometry, including B\'ezout theorem, a discussion of birationnal equivalence,  and  the  quartic associated to the double equations of degree 2. He gives some examples of applications of elementary algebraic geometry to diophantine analysis avoiding the megamachines often needed in these kind of results.

In Number Theory, 
he worked on Waring's problem.
He wrote an article {\sl Une pr\'esentation ad\'elique de la s\'erie singuli\`ere et du probl\`eme de Waring
[An adelic presentation of the singular series and the Waring
problem]} \cite{rod-lac5} to
 consider Waring's problem and the corresponding singular series in the
adelic setting from ideas of Igusa, Ono and Weil.
The local singular series in the place $p$ is expressed in the form of counting points.
Later it is this expression that applies in the demonstration of some examples of Yuri Manin's conjectures about the asymptotic behaviour of the number of points of bounded height on an algebraic variety defined over  a number field when the height goes to infinity by  some people who tackled that subject such as Antoine Chambert-Loir, Yuri Tschinkel, R\'egis de la Bret\`eche,  Tim Browning, Emmanuel Peyre.

He worked also on real quadratic fields.
In the article {\sl On real quadratic fields}~\cite{rod-lac7},
 he studies 
the caliber of a real quadratic field $Q(\sqrt d)$ (the number $k(d)$
of reduced primitive binary quadratic forms of discriminant $d$).
A theorem of Siegel says
 there are but a finite number of real quadratic fields with
a given caliber but it is not effective.
Gilles shows that under GRH, the only real quadratic fields with
caliber one are the seven fields $Q(\sqrt t)$ with $t = 2, 5,13,29, 53, 173, 293$. 
Also he shows that under GRH, the
only principal real quadratic fields with discriminant $d = r^2 +1\equiv 1 \pmod 4$ are the six fields $Q(\sqrt t)$
with $t = 5, 17, 37, 101, 197, 677$.
For that he proved a formula involving the class number $h(d)$ and the fundamental unit $e_0$:
$h(d) \log e_0 (d) < c \cdot k(d)$ where $c <4.230$.
In 2010, this has been proved and generalised by 
Byungheup Jun and Jungyun Lee, in the article {\sl Caliber number of real quadratic fields} \cite{JunLee}, where they find, without GRH,  the fields with caliber 1 and also, 
the fields of caliber 2 and discriminant $\not\equiv5 \pmod 8$.

We did not see each other so often anymore  when he was appointed professor in Nice in 1981. 
But an article of Tsfasman and Vl\v{a}du\c{t} dealing with application of modular functions to the theory of Goppa codes led Gilles to be interested in Information Theory.
 In 1985, he  made a presentation at the prestigious Bourbaki seminar on {\sl the Geometric  Goppa codes} \cite{rod-lac8} where he explained these works. I was there, especially since at that time, I wanted to change my field of research from the Theory of Group Representations to the Theory of Codes.

It was also in 1985 that Gilles took up the position of \emph{Directeur de Recherches CNRS} and moved from Nice to Luminy. 
He became the  second director of CIRM, from 1986 to 1991. The CIRM was inaugurated in 1991, and Gilles knew it  already, as  vice-president of SMF, i.e., {\sl Soci\'et\'e Math\'ematique de France} (French Mathematical Society), during 1982--1984, and as the  president of the scientific council of CIRM. 
His stewardship was around the beginning of CIRM. 
 With the technical team, he had a pioneering spirit.
During that period,
there was a rise in scientific activities: 
from 20 symposia in 1986, it reached  37 in 1991.

Already, at that time, a certain number of symposia foreshadowed in part the major currents of mathematical research.
There were some big conferences, such as the 300-person symposium organized by Mohammed Mebkhout, or the {\sl Journ\'ees Arithm\'etiques} 
in 1989, with a similar number of participants, organised by Gilles.
The highlight of this period was the construction of the new library. The former was housed in the premises of the University library,  it was open at the same hours  and it was far from CIRM. Obtaining the necessary funding to build new library for  CIRM was not easy, but finally he was able to launch an architectural competition.
Of course there was some complication during the works. For instance
he  wanted to build a beautiful amphitheater in the new building, but while  digging, they fell on the rock of a creek (that's why the conference room at that time was horizontal).
The library was finally inaugurated in 1991, at the time of the 10th anniversary of CIRM. It had been a really short time (5 years)
since the granting of the credits until the inauguration.
Gilles  showed his ecological feeling and ended this period with a long-term wish: CIRM has 4 hectares of a forest that is quite exceptional, since it is a northern forest in the middle of the southern forest of creeks. The flora is well supplied, but he complained that the forest areas are not restored at all. 
He wished those who follow him to achieve preservation and enrichment 
 of this precious ecological heritage.

In 1987, while Gilles was the head of CIRM, he organised with Jacques Wolfmann from the University of Toulon, a conference named Algebraic Geometry and Coding Theory (AGCT). At that time, one did not imagine that this conference would still exist now under a slightly modified name because of the renewal of the subjects: Arithmetic, Geometry, Cryptography, and Coding Theory (AGC$^2$T).

He did several works concerning popularization of Mathematics.
He wrote several articles in the general public journal {\sl La Recherche} (Les codes correcteurs d'erreurs) or for the {\sl Encyclopedia Universalis} (Qu'est ce qu'une \'equation).

In the article
{\sl The
weights of the orthogonals of the extended quadratic binary Goppa
codes} \cite{rod-lac9} with Jacques Wolfmann
they started from results on elliptic curves and Kloosterman sums over the finite field $\mathbb{F}_{2^t}$, and they determined the weights of the orthogonals of some binary linear codes.
The function $f (x) = x^{-1}$ on $\mathbb{F}_{2^t}$  resists to differential cryptanalysis and it  is used in the ``Rijndael" protocol retained in the call for proposals for the Advanced Encryption Standard (AES) launched by the National Institute of Standards and Technology (NIST).
They could compute
the Walsh spectrum of that function  from their result giving a relation between the number of points of ordinary elliptic curve and  Kloosterman sums.
This article has become the most quoted article of his.

A little before 1990, Gilles created a working group on the Theory of Information that he set up at CIRM, 
that included 
Guy Chass\'e, Yves Driencourt, Dominique Le Brigand, Jean-Francis Michon, Marc Perret, Robert Rolland, and myself. This group was very active and organised a meeting every month.

At that time, he was also working with G\'erard Rauzy and Jean-Yves Girard and with the head of CNRS to create a CNRS proper laboratory dedicated to Discrete Mathematics. With CIRM, CNRS wanted to create a ``mathematics  campus" in Luminy. The LMD (Discrete Mathematics Laboratory) was created in Marseille in October 1992 with a team ATI (Arithmetic and Theory of Information) that he led and that I joined  naturally.

Gilles was
generous, always willing to help
others.
He helped me a lot to make my moving easier, because I was disabled at the time. 
In order that I get used to giving talks to students he left me a course of DEA/Master2 to teach in his place. 

He became the Director of IML (Mathematics Institute of Luminy) from 2006 to 2011, 
and 
then gave me the place of Director. 
At that time, I saw the investment in time and energy that 
this appointment required, and wondered how he had managed 
to occupy 
this position twice while he kept on doing mathematics.

He was optimistic as shows the last letter I received from him, one month before he died.

\begin{quotation}
{\it ``Pour ma sant\'e, \c ca va bien, je suis normal, j'ai un excellent traitement,
[...]
Le r\'esultat est que je deviens casanier et paresseux.
Donc je fais des maths \`a la maison, sur les statistiques de traces de repr\'esentations, ce qui est une activit\'e qui se fait assis.

J'arrive quand m\^eme \`a sortir dans le quartier, et j'ai pu aller \`a Besan\c con pr\'esider la th\`ese de Philippe Lebacque.

Il faudrait tout de m\^eme qu'on se voie un de ces jours. Quel est ton jour? Le jeudi?"
}\footnote{\it ``For my health, I'm fine, I'm normal, I have an excellent treatment,
[...]
The result is that I become a homebody and lazy.
So I do maths at home, on statistics of traces of representations, which is an activity that is done sitting.

I can still go out in the neighborhood, and I went to Besan\c{c}on to preside over the thesis of Philippe Lebacque.

We should still see each other one of these days. What is your \hbox{day?} Thursday?"}
\end{quotation}

Besides mathematics, Gilles's interests included philosophy, literature, art, history, and so on.
We miss Gilles a lot. He breathed tranquility. I retain his vast culture and his humanism. He leaves a great void.
%
%

\medskip

\section{Gilles Lachaud: friend and mathematician}
\medskip
\centerline{Michael A. Tsfasman}
\medskip

Algebraic geometry codes were discovered by Goppa at the beginning of 1980s and 
became widely known to mathematicians after 1981 when our preprint of \cite{TsVlZ}---with the proof that 
there exist codes over the Gilbert--Varshamov bound---was smuggled through the iron curtain. In 
1985, Gilles Lachaud gave a talk on the subject \cite{rod-lac8} 
at the Bourbaki seminar. Most likely he was attracted by the subject because our proof was based on modular curves, and his basic education was in modular and automorphic forms.  In 1987, he organised the first AGCT meeting at CIRM and invited me to come. Strong mutual dislike between me and the Soviet regime naturally led to the refusal of the leave permit needed at the time to cross the border. Two years later---the regime being much weakened by perestroika---I came for AGCT-2, and we saw each other for the first time.
Gilles met me at the entrance to CIRM, and we became close friends ever after.

I was impressed by the scope of his interests with quite an intersection with mine: literature, philosophy, history, religion, art, and science, not to speak about mathematics. I appreciated it even more when I learned French enough to switch to it in our communication.

Gilles was the director of CIRM for 6 years and his influence is felt there even now. He built the library that should bear his name. And he reorganized CIRM in such a way as to make it maximally friendly for the mathematical community in everything from perfect blackboards to wonderful cuisine. 
Then he organised the LMD (Discrete Mathematics Laboratory) that later became the Institute of Mathematics of Luminy.

It is because of him that among the countries to work outside Russia I have chosen France.

Now let me describe the part of his legacy that is close to my interests. I mean mostly geometry over finite fields and a little bit coding theory.

Gilles's work on codes and curves over finite fields started with questions where he could use his expertise in exponential sums, Kloostermann sums, Eisenstein sums, and so on \cite{Eisenstein, LW, CU, CU2, AS, Me}.

He also did some nice work, joint with J. Stern, on purely coding theory questions \cite{LS1, LS2}. In the first of these papers they prove that there are families of asymptotically good codes over $\mathbb F_q$ with
polynomial complexity of construction whose relative weights are as close to $(q - 1)/q$
as we want. In the second one they prove a nice asymptotic lower bound for polynomially constructed spherical codes and for the (polynomial) kissing number. 

How many $\mathbb F_q$-points can there be on an abelian variety? On the Jacobian of a curve? Since $h_X=\vert J_X (\mathbb F_q)\vert=P_X(1)$ the obvious bound is 
$$
(\sqrt{q}-1)^{2g}\le h \le (\sqrt{q}+1)^{2g}.
$$
In his paper 
with Martin--Deschamps \cite{LMD}, they prove much better bounds:

\medskip

{\sc Theorem.} 
\emph{For a smooth geometrically irreducible algebraic curve $X$ of genus $g$, the number $h_X=\vert J_X (\mathbb F_q)\vert$ satisfies the lower bounds
$$
h_X\ge q^{g-1}\hskip 0.1cm
 \frac{(q-1)^2}{(q+1)(g+1)}
$$
and 
$$
h_X\ge (\sqrt{q}-1)^2 \hskip 0.2cm\frac{q^{g-1}-1}{g}\hskip 0.2cm\frac{\vert X (\mathbb F_q)\vert+q-1}{q-1}. 
$$
If $g>\frac{\sqrt{q}}{2}$ and if $X (\mathbb F_q) $ is nonempty, then
$$
h_X\ge {(q^g-1)}\hskip 0.1cm\frac{q-1}{q+g+gq}. 
$$
}

Much later, in joint papers with Aubry and Haloui \cite{AHL, AHL'}, they ameliorated these bounds. In the same paper  they proved an analogue of Serre's bound for an arbitrary abelian variety.

His paper with  S. Vl\v{a}du\c{t} \cite{LV} treats the function field analogue of the famous Gauss problem asking whether the number of real quadratic fields with class number 1 is infinite, or not. Its weak version asks whether the number of principal algebraic number fields (i.e., with class number 1) is infinite. Both are wide open. In the function field case we have two types of their analogues, the obvious one asks just for infinity of such fields and is easy to solve in the affirmative, the more reasonable one asks for infinity of such fields when the constant field $\mathbb F_q$ is fixed.

Here is an illustrious specimen of their results:

\medskip

{\sc Theorem.} 
\emph{
Assume $q=4, 9, 25, 49$, or $169$. Let $P\in \mathbb {P}^1(\mathbb{F}_q)$. Then
there are infinitely many extensions $(K, S)$ of the pair $(\mathbb{F}_q(T), \{P\})$ such
that:}

(1) \emph{The place $P$ completely splits in $K$.}

(2) \emph{The field $K$ is a Galois extension of $\mathbb{F}_q(T)$.}

(3) \emph{The ring $A_S$ is principal.}

\medskip

Now let me discuss in more details our only joint paper {\sl Formules explicites pour le nombre de points des vari\'et\'es sur un corps fini [Explicit formulas for the number of points of varieties over a finite field]} \cite{LTs}. In mid-nineties, I realised that though there are plenty of mathematicians working on curves over finite fields and their number of points, actually no-one had ever studied the number of points on surfaces. After a while I had decided to write down my first results on surfaces and multi-dimensional varieties. It was the moment of my intensive attempts to study French (Gilles' influence, of course). So I started to write my only French language paper \cite{Ts}. Quite naturally, the French of it was awful. So I asked Gilles to correct it. He did it perfectly and little by little we started to think on the subject together.

For curves the asymptotic problem is that of the behaviour of the ratio $\vert X (\mathbb F_q)\vert/g_X$. 
For multidimensional varieties I suggested to divide by the sum of Betti numbers. The problem becomes mutidimensional, say for surfaces we have two parameters instead of one, $x=b_1/(b_1+b_2)$  and $y=N/(b_1+b_2)$, here $N=\vert X (\mathbb F_q)\vert$. Therefore, the answer is not a number as in the case of curves (where, at least when $q$ is a square,  it is equal to $\sqrt{q}-1$), but a curve on the $(x,y)$-plane.

One of the results of that paper of mine was the generalization of the Drinfeld--Vl\v{a}du\c{t} asymptotic upper bound to the multi-dimensional case. I did that using the elementary approach of Drinfeld--Vl\v{a}du\c{t}. It was clear for both of us that the explicit formulae approach could be better. 

The process of our work showed the difference between the Russian and French mathematical schools, as well as Gilles' style of doing mathematics. I would give a vague and often incorrect idea and try to explain it to him. Next day he would come with a couple of pages: precise definitions and clear results. I would read and criticize it, and then in the morning a better version was given to me. 

After you write down explicit formulae, to get bounds you need to choose good test functions. These test functions should be doubly positive, meaning that both their values and their Fourier coefficients are positive. Gilles noticed that these functions happen to be related with (doubly) positive kernels used in sophisticated analysis starting from 19th century. After handing the previous night results to me to read, he would go to CIRM library to look through classical books to find new and new ones for us to use.

Classics mainly considered doubly positive kernels (densities of the corresponding measures) such as, for example, Fej\'er, Johnson, and de la Vall\'ee-Poussin kernels, and the Jacobi kernel.

There is also a more recent source of doubly positive sequences. Namely, J.~Oesterl\'e, who studied for the case of curves a more subtle question on the minimum genus of a curve having a prescribed number $N$ of points, used explicit formulas and substituted into them doubly positive sequences optimal for precisely this problem, which involves two free parameters, $N$ and $g$. However, in the asymptotic problem, two parameters appear not for curves but for surfaces, so Oesterl\'e sequences yield good bounds for the case of surfaces too.

Each of the above-mentioned doubly positive kernels or sequences gives a specific bound on the number of points. There are
countably many such bounds, their graphs intersect each other, and the best known bound is given by their envelope curve.


\begin{figure}[tp]
\centering\vskip5pt \unitlength=.001\textwidth\small
\begin{picture}(0,0)
\put(12,-5){0} \put(12,86){1} \put(12,177){2} \put(12,268){3} \put(12,359){4}
\put(106,-22){0{,}2} \put(195,-22){0{,}4} \put(283,-22){0{,}6} \put(372,-22){0{,}8}
\put(471,-22){1}
\end{picture}
\hspace*{\fill}\includegraphics[width=.45\textwidth]{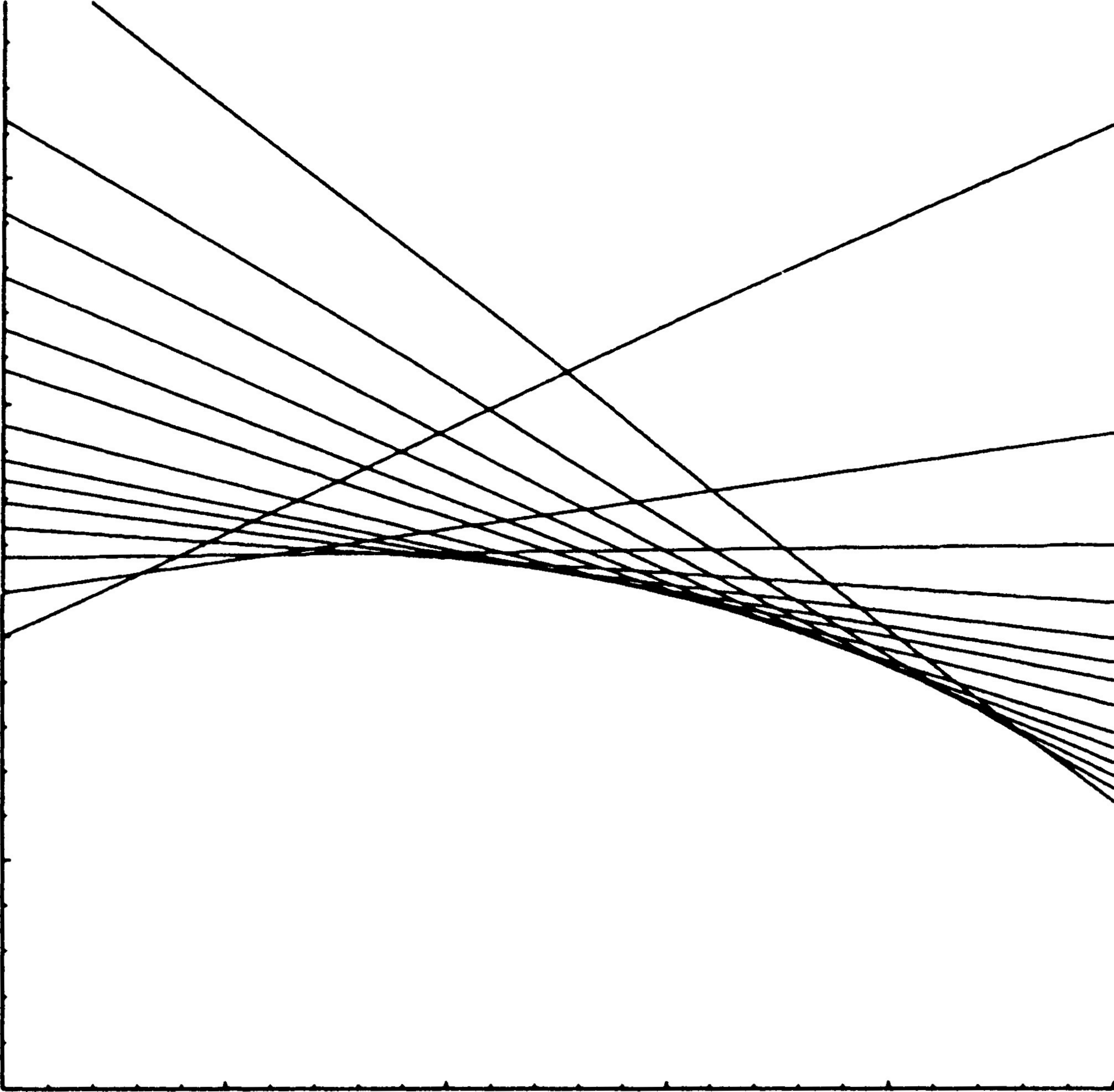}\hfill
\begin{picture}(0,0)
\put(12,-5){0} \put(12,86){1} \put(12,177){2} \put(12,268){3} \put(12,359){4}
\put(105,-22){0{,}2} \put(193,-22){0{,}4} \put(281,-22){0{,}6} \put(370,-22){0{,}8}
\put(468,-22){1} \put(473,395){$W$} \put(473,285){$S$} \put(473,125){$L$}
\put(473,40){$M$}
\end{picture}
\hfill\includegraphics[width=.45\textwidth]{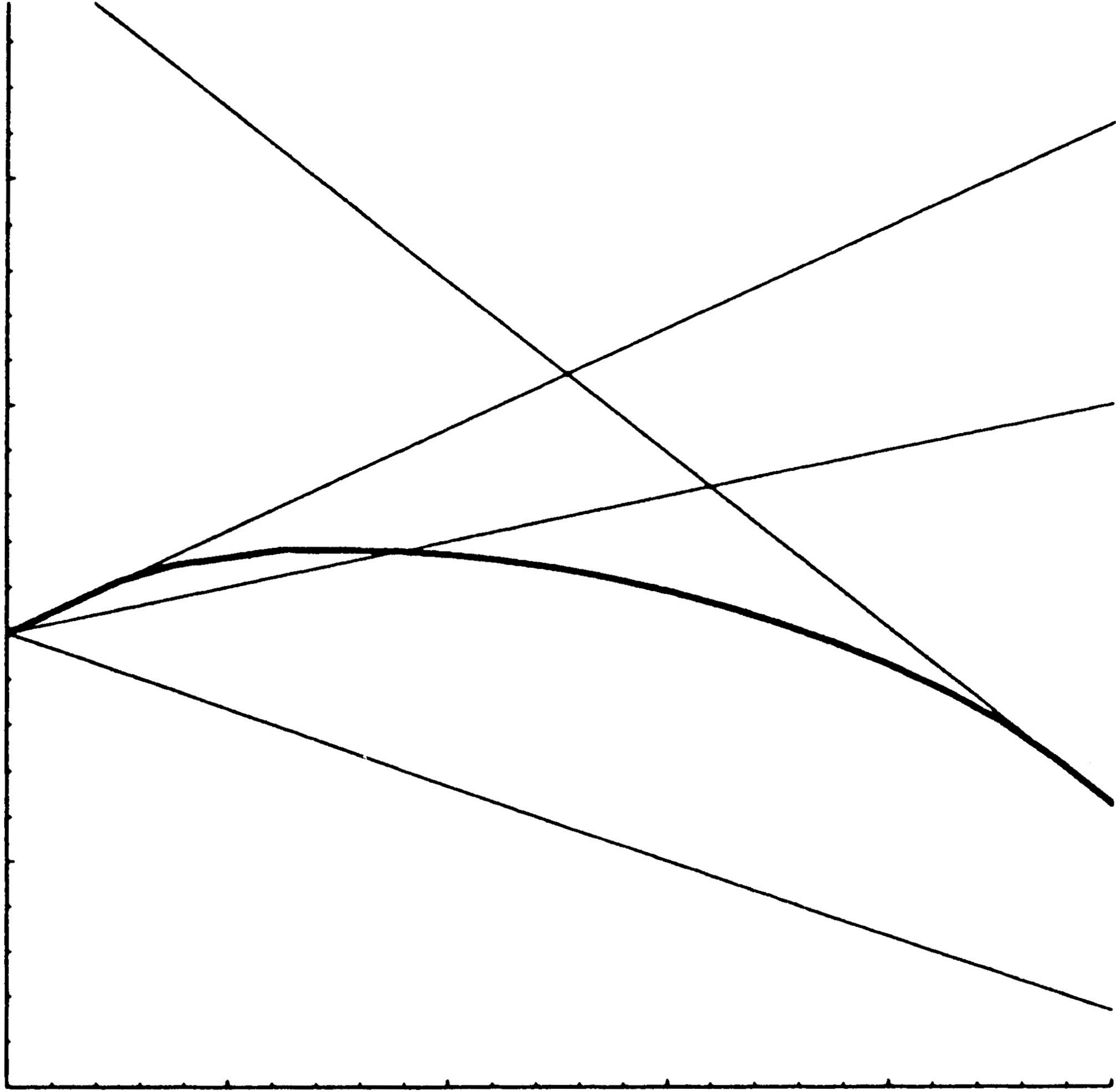}\hskip2pt\hspace*{\fill}\strut\\
\vskip20pt \hspace*{\fill}Fig.~1.\hfill\hfill Fig.~2.\hspace*{\fill}
\end{figure}

The result (one of, the paper being large enough) is depicted on the following figures for the case of $q=2$ (on these plots, the abscissa axis corresponds to the parameter $\frac{b_1}{b_1+b_2}$, and the ordinate axis, to the parameter
$\frac{N}{b_1+b_2}$): Fig.~1 shows bounds that correspond to the above-described doubly positive objects and their envelope; Fig.~2, besides the envelope, shows the Weil bound~$W$, Serre bound~$S$, generalized elementary bound~$L$, and the existence bound~$M$ from \cite{Ts}.

\medskip

\section{A Tribute to Gilles Lachaud}
\medskip
\centerline{Sudhir R. Ghorpade}

\medskip

During August 
1993, there was a 3-week CIMPA\footnote{CIMPA =  
\emph{Centre International de Math\'{e}matiques Pures et Appliqu\'{e}es} = \emph{International Centre of Pure and Applied Mathematics}. It is an organization based in France that conducts instructional schools on mathematics in ``developing countries". For more information, see: \href{https://www.cimpa.info/en}{https://www.cimpa.info/en}} School on Discrete Mathematics at Lanzhou, China. The speakers included Gilles Lachaud, Maurice Mignotte, Robert Rolland, and Wu Wen-Tsun.  I was one of the participants in this school, and it was here that I first met Gilles Lachaud. As it turned out, 
that was the beginning of a long and cherished mathematical association and close friendship. 
In what follows, we will discuss some of Gilles' work on topics such as counting points of algebraic varieties over finite fields, and linear error correcting codes. 

In the previous two sections, the Bourbaki seminar \cite{rod-lac8}  and some of the earlier works of Lachaud on coding theory, including his much cited paper with Wolfman \cite{rod-lac9}, have been mentioned. I will begin with an innocuous little paper on projective Reed-Muller codes that first appeared in the Cachan proceedings \cite{RM}. 
Recall that given any positive integers $d,m$,  the (generalized) Reed-Muller code $\RM_q(d,m)$ of order $d$ and length $n=q^m$ can be defined as the image of the evaluation map 
$$
\Ev : \Fq[x_1, \dots , x_m]_{\le d} \to \Fq^n \quad \text{given by} \quad \Ev (f) = (f(P_1), \dots , f(P_n)), 
$$
where $ \Fq[x_1, \dots , x_m]_{\le d}$ denotes the $\Fq$-vector space of all polynomials of (total) degree $\le d$ in $m$ variables $x_1, \dots , x_m$ with coefficients in the finite field $\Fq$, and where $P_1, \dots , P_n$ is an ordered listing of distinct points in $\Aff^m(\Fq)$. This has been widely studied in coding theory, and we cite in particular, the work of 
Delsarte,  Goethals, and MacWilliams \cite{DGM}. Gilles proposes a natural generalization of this where the affine space $\Aff^m$ is replaced by the projective space $\PP^m$ and quite naturally, $\Fq[x_1, \dots , x_m]_{\le d}$ is replaced by the space $\Fq[x_0, x_1, \dots , x_m]_{d}$ of homogeneous polynomials of degree $d$ in $m+1$ variables. Further, he takes fixed representatives in $\Fq^{m+1}$ of the points of $\PP^m(\Fq)$ (for instance, those where the first nonzero coordinate is $1$). One can then define a similar evaluation map and its image is the projective Reed-Muller code $\PRM_q(d,m)$ of order $d$. It has length
\begin{equation}
\label{eq:pm}
p_m:= | \PP^m(\Fq)| = q^m + q^{m-1} + \dots + q +1.
\end{equation}
Moreover, if $d\le q$, then the evaluation map is injective and so the dimension of $\PRM_q(d,m)$ is ${{m+d}\choose{d}}$. In \cite{RM}, Gilles gives a lower bound on the minimum distance of $\PRM_q(d,m)$, and computes the exact value  in the special cases $d=1$ (where the Plotkin bound is attained) and $d=2$. The latter uses classical facts about zeros of quadrics in finite projective spaces. For the general case, one needs tight bounds on the maximum number of $\Fq$-rational points on a hypersurface in $\PP^m$ or in other words, the maximum number, say $N_d$, of zeros in $\PP^m(\Fq)$ that a homogeneous polynomial of degree $d$ in $m+1$ variables with coefficients in $\Fq$ can have. Tsfasman conjectured that if $d \le q+1$, then
$
N_d = dq^{m-1} + p_{m-2},
$
where $p_j$ is as in \eqref{eq:pm} if $j$ is a nonnegative integer and $p_j:=0$ if $j <0$. Tsfasman's conjecture was quickly proved in the affirmative by J.-P. Serre, who communicated it to Tsfasman in a letter which was published some two years later \cite{Se}. Gilles was then in a position to complete the results of \cite{RM} and give a formula for the minimum distance in the general case in \cite{RM2}, where he argues that from coding theoretic viewpoint, projective Reed-Muller codes are somewhat better than (affine, or generalized) Reed-Muller codes, and as such, they deserve to be studied further. A beginning was already made by Sorensen \cite{So} whose paper, published in the same year as \cite{Se}, and gives, in effect, an alternative proof of Tsfasman's conjecture and additional results of particular interest to coding theory. 

The introduction of projective Reed-Muller codes by Gilles and the question about the minimum distance of these codes, which is essentially equivalent to Tsfasman's conjecture mentioned above, has had many offshoots, and we will briefly describe some of them. 

Tsfasman's conjecture is about projective hypersurfaces defined over $\Fq$ of degree $d \le q+1$. Now a hypersurface is the simplest example of a complete intersection in $\PP^m$ of high dimension (vis-\`a-vis $m$). 
Keeping this in view, Gilles formulated the following conjecture (which appeared in print much later, in our joint paper \cite{GL02}).  

\medskip

{\sc Conjecture.}
\emph{If $X$ is a complete intersection in $\PP^m$ defined over $\Fq$ of dimension $n \ge {m}/{2}$ and degree 
${d } \le q + 1$, then} 
$$
 |X(\Fq)| \le  
{d } p_n - ({d }-1) p_{2n - m} =  {d }(p_n - p_{2n - m}) + p_{2n - m}.
$$

Observe that if $X$ is a hypersurface (so that $n = m-1$), then this is {Tsfasman's Conjecture}. 
On the other hand, Tsfasman together with Boguslavsky, generalized his conjecture about hypersurfaces of a given degree $d \le q+1$ to intersection, say $S_{d,r}$ of $r$ hypersurfaces, each defined by a  homogeneous polynomial of degree $d < q$, such that these $r$ polynomials are linearly independent. Tsfasman-Boguslavsky Conjecture proposes an intricate, but explicit, formula $T_r(d,m)$ for the maximum number of $\Fq$-rational points of algebraic sets such as $S_{d,r}$. Again, if $r=1$, then this reduces to Tsfasman's Conjecture about hypersurfaces. 

In 2015-2017, that is, roughly about two decades after these conjecture were made, two interesting developments took place almost simultaneously. On the one hand, Alain Couvreur \cite{Co16} showed that  Lachaud's Conjecture, and in fact,  a more general result, holds in the affirmative. On the other hand, Mrinmoy Datta and I showed that the 
Tsfasman-Boguslavsky Conjecture holds if $r\le m+1$, but it  is false, in general (cf. \cite{DG15}, \cite{DG17}). I think Gilles was rather happy and proud that his guess turned out to be correct. Shortly before these developments took place, Gilles himself returned to the topic and together with Robert Rolland, he wrote a nice article \cite{LR15} giving various bounds for the number of $\Fq$-rational points of affine as well as projective varieties over finite fields. 
In AG${\rm C}^2$T-2015, Lachaud, Couvreur and I gave talks in a single session to explain these developments, and I remember how Gilles provided a nice perspective on the topic in the first talk of that session. I might mention in passing that a refined version of Tsfasman-Boguslavsky Conjecture has been proposed in \cite{DG17} and \cite{BDG}.  The newer conjectures remain open in general although they are known to be valid in many cases (cf. \cite{BDG18} and \cite{BDG}). 

Another interesting offshoot of the work of Gilles on projective Reed-Muller code and the inequality Serre proved to answer Tsfasman's Conjecture is a development where Gilles himself was involved. It is the subject matter of his  paper \cite{IPAM}, which was published barely six months before he passed away. Perhaps it is worthwhile to mention a background to this. In February 2016, 
Everett Howe, Kristin Lauter, and Judy Walker organized a workshop on Algebraic Geometry for Coding Theory and Cryptography at the Institute for Pure and Applied Mathematics (IPAM) on the campus of the University of California, Los Angeles. They said that this will not be a standard ``five talks per day with a break for lunch" type of conference, but instead an activity focused on research projects that will be started during the workshop with a hope that it will result in ongoing collaboration among participants. Gilles and I were asked to be co-leaders of a working group. Both of us were apprehensive about this experiment. Nonetheless, we put together a group of six participants including Yves Aubry, Wouter Castryck, Mike O'Sullivan, and Samrith Ram. The broad topic for our group was ``Number of points of algebraic sets over finite fields".  After some introductory talks, we discussed every day of the week (and sometimes late in the evenings as well) and bounced ideas off each other. As a result of this brainstorming, and continued correspondence thereafter, we wrote an article concerning an analogue of Serre's inequality for hypersurfaces in weighted projective space, and weighted projective Reed-Muller codes. Gilles was a very active participant in this, and wrote almost single handedly a compendium on the (geometry of) weighted projective spaces as a handy manual for us, which formed an appendix to our article \cite{IPAM}. Although there are many partial results in this article, the problem of determining the maximum number of $\Fq$-rational points on a hypersurface in weighted projective spaces is open, in general, and a conjecture stated here could be a pointer for further research on the topic. Before ending the discussion of projective Reed-Muller codes introudced by Gilles, I would like to metion an article of well-known finite geometers and group theorists Bill Kantor and Ernie Shult 
\cite{KantorShult}, which was published in 2013. Here they remark that their first main theorem can be viewed as a statement about a certain code $C$ having a check matrix whose columns consist of one nonzero vector in each Veronesean point, Then they write: \emph{We have not been able to find any reference to this code in the literature. It is probably worth studying, at least from a geometric perspective.} In fact, the code $C$ that they talk about is nothing but the dual of the projective Reed-Muller code studied by Gilles about 
25 years ago!

I would now like to go back to my first meeting with Gilles at the 1993 CIMPA workshop in China, which was also the first time I learned about error correcting codes from the lecture courses of Robert Rolland and Gilles Lachaud. In the course of lectures or possibly, during the discussions over meals, I learned about the seemingly intractable open problem of determining an explicit formula for the number, say $\gamma (q; k,n)$, of $q$-ary MDS codes of length $n$ and dimension $k$. The answer was known for $k\le 2$ and also for $k=3$ and $n\le 8$. I could see that the problem was equivalent to determining the number of $\Fq$-rational points of the open subset $U_{k,n}$ of the Grassmannian $G_{k,n}$ consisting of points of $G_{k,n}$ for which \emph{all} the Pl\"ucker coordinates are nonzero. Based on this observation, I could make some observations, which interested Gilles and he encouraged me to think further. The activity picked up when Gilles arranged for me to get a fellowship from CNRS to visit him Marseille for a month in June 1996. This was the first of my numerous visits to Marseille in subsequent years. Using the connection with Grassmannians together with the work of Nogin \cite{No96} on Grassmann codes, and various combinatorial inputs, we could obtain explicit lower and upper bounds for $\gamma (q; k,n)$ and an asymptotic formula in the general case. Further, we could use some of the auxiliary results together with Grothendieck-Lefschetz trace formula to derive some geometric properties of certain linear sections of Grassmannians. These results are given in our article {\sl Hyperplane sections of Grassmannians and the number of MDS linear codes}  \cite{GL01}. Originally, there was to be a remark in this paper indicating how some  results of Nogin on higher weights of Grassmann codes could also be derived as a consequence of the combinatorial results we had in preparation for our results on $\gamma (q; k,n)$. But then there was to be a conference in Guanajuato, Mexico in 1998, which Gilles encouraged me to attend. We decided to expand this remark and used the opportunity to introduce what seemed a natural generalization of Grassmann codes, called Schubert codes. This paper {\sl Higher weights of Grassmann codes} \cite{GL00} was published in the proceedings of the Guanajuato conference and it contained a conjecture about the minimum distance of Schubert codes. As it turned out, \cite{GL00} led to a lot of further research and it has many more citations than \cite{GL01}. The moral of the story for me (and I think Gilles agreed with that) was: A good problem is sometimes more valuable than a good theorem! 

In the course of our work on MDS codes, we thought it would be nice to have concrete estimates for the number of $\Fq$-rational points of arbitrary algebraic varieties, especially when the variety is not necessarily smooth. The best general estimate we had at our disposal was the following inequality proved by Lang and Weil in 1954:

\medskip

{\sc Lang-Weil Inequality.} 
\emph{Let $X$ be an (absolutely) irreducible projective
variety} in $\PP^N$ defined over $\Fq$. If $X$ has dimension $n$ and degree $d$, then
$$
\large| {|X(\Fq)| - p_n }\large| \le (d-1)(d-2) q^{n - (1/2)} + C q^{n-1} ,
$$
\emph{where $C$ is a constant depending only on $N$, $n$, and $d$.
}

\medskip

When the variety is smooth, and better still, a complete intersection, much sharper estimates are available such as the one obtained by Deligne in 1973 as a consequence of his seminal work on Weil conjectures. 

\medskip

{\sc Deligne's Inequality for Smooth Complete Intersections.}  
\emph{If $X$ is a  {nonsingular complete
intersection} in $\PP^N$ over $\Fq$ 
of dimension $n = N-r$, then  
$$
\large|{ |{X(\Fq)}| - p_{n} } \large|  \leq b'_n \, q^{n/2}.
$$
Here $b'_n = b_n - \epsilon_n$ is the {primitive $n$th 
Betti number} of $X$ (where $\epsilon_n = 1$ if $n$ is even and $\epsilon_n=0$ if $n$ is odd), and $p_n:=|\PP^n(\Fq)| = q^n + q^{n-1} + \cdots + q +1$. 
}

We~remark that if $X$ has multidegree 
$\mathbf{d} = (d_1, \dots , d_r)$, 
then a formula of Hirzebruch shows that $b'_n = b'_{n}(N, \mathbf{d})$ equals 
$$
\textcolor{black}{(- 1)^{n + 1} (n + 1) +
\sum_{c = r}^{N} (- 1)^{N+c} \binom{N + 1}{c + 1}
\mathop{\sum_{\nu_1+ \cdots + \nu_r = c}}_{\nu_i\ge 1 \; \forall i}    d_1^{\nu_1} \cdots  d_r^{\nu_r} }. 
$$ 

From a practical viewpoint, one likes to have some control on the constant $C$ appearing in Lang-Weil inequality and have a Deligne-like inequality for complete intersections (in particular, hypersurfaces) that are not necessarily smooth. We had an idea how this might be achieved, and we worked for years together to prove the following result (\cite[Theorem 6.1]{GL02}), which is one of the main results of our paper {\sl \'Etale cohomology, Lefschetz theorems and number of points of singular varieties over finite fields} \cite{GL02}, published in a special volume of \emph{Moscow Mathematical Journal} dedicated to Yuri I. Manin on the occasion of his 65th birthday. 

\medskip

 {\sc Deligne-type inequality for arbitrary complete intersections.}
\emph{  
Let $X$ be an irreducible complete intersection of dimension $n$ in
$\PP^N_{\Fq}$, defined by $r = N - n$ equations, with multidegree
$\mathbf{d} = (d_{1}, \dots, d_{r})$, and let $s$ be 
the dimension of the singular locus of  $X$. 
Then
$$
| {|{X(\Fq)}| - \, p_{n}}| \leq
b'_{n - s - 1}(N - s - 1,\mathbf{d}) \, q^{(n + s + 1)/2} +
C_{s}(X) q^{(n + s)/2},
$$
where $C_{s}(X)$ is a constant independent of $q$. If $X$ is nonsingular, then
$C_{- 1}(X) = 0$. If $s \geq 0$, then}
$$
C_{s}(X) \leq 9 \times 2^{r} \times (r \delta + 3)^{N + 1}
\quad {\it where } \quad \delta = \max\{d_{1}, \dots, d_{r}\}. 
$$

\medskip

For normal complete intersections, this may be viewed
as a common refinement of Deligne's inequality and the Lang-Weil inequality.  Corollaries include previous 
results of {Aubry and Perret} (1996), {Shparlinski\u{\i} and Skorobogatov} (1990), as well as 
{Hooley and Katz} (1991). 

We also used the power of Weil conjectures and some estimates of Katz (2001) to arrive at a version of Lang-Weil inequality with an explicit bound on the constant $C$ appearing therein (\cite[Theorem 11.1]{GL02}). A precise statement is given below. Here, given any $m$-tuple $\mathbf{d} = (d_{1}, \dots, d_{m})$ of positive integers, we say that an affine (resp. projective) variety in 
$\Aff^N_{\Fq}$ (resp. $\PP^N_{\Fq}$) is of \emph{type} $(m, N, \mathbf{d})$ if it can be defined by the vanishing of $m$ polynomials in $N$ variables (resp. $m$ homogeneous polynomials in $N+1$ variables) with coefficients in $\Fq$. 

\medskip

{\sc Effective Lang-Weil inequality.} 
\emph{Suppose $X$ is a projective variety in $\PP^N_{\Fq}$ or an affine
variety in $\Aff^N_{\Fq}$ defined over $\Fq$. Let $n = \dim X$
and $d = \deg X$. Then
$$
|{|{X(\Fq)}| - p_{n}}| \leq (d - 1)(d - 2) q^{n - (1/2)} +
C_{+}(\bar{X}) \, q^{n - 1},
$$
where $C_{+}(\bar{X})$ is independent of $q$. Moreover, if $X$ is of type 
$(m, N, \mathbf{d})$, with $\mathbf{d} = (d_{1}, \dots, d_{m})$, and if 
$\delta = \max\{d_{1}, \dots, d_{m}\}$, then}
$$
C_{+}(\bar{X}) \leq \begin{cases}
9 \times 2^{m} \times (m \delta + 3)^{N + 1} & {\it if} \ X \ {\it is \ projective} \\
6 \times 2^{m} \times (m \delta + 3)^{N + 1} & {\it if} \ X \ {\it is \ affine.} \\
\end{cases} 
$$

\medskip

As a corollary, one obtains an analogue of a result of {Schmidt} (1974) on a lower bound
for the number of points of irreducible hypersurfaces over $\Fq$. 

Finally, we showed that certain \emph{conjectural statements} of Lang and Weil (1954), 
relating the Weil zeta function of a smooth projective variety $X$ defined over $\Fq$ and the ``characteristic polynomial'' of its Picard variety, hold true provided one uses the ``correct'' Picard variety.  For more details, one may refer to our paper \cite{GL02} or the expository article {\sl Number of solutions of equations over finite fields and a conjecture of Lang and Weil} \cite{GL02b} that appeared in the proceedings of an international conference on 
Number Theory and Discrete Mathematics, held at Chandigarh, India, which Gilles attended. 

Although it took a long time to complete this work, it was quite satisfying in the end. What was particularly gratifying was to see many seemingly disparate topics to which our results found applications (references to some of such applications as  well as to extensions and generalizations of our results, especially by Cafure and Matera, can be found in the corrigenda and addenda  \cite{GL02c} to \cite{GL02}). For me, it was a great learning experience to work with Gilles on this topic and to witness firsthand his tenacity and technical skills. At the beginning of this paper, there is a \emph{Sanskrit} verse from \emph{Rg Veda} that can roughly be translated as ``Their cord was extended across". Now one might be tempted to think that this was put in there by the Indian coauthor. But the fact is it was Gilles who insisted that we insert it. To him, it conveyed somehow that we are merely extending the key ideas of our illustrious 
predecessors. Also, he felt it was pertinent because we frequently use Bertini-type theorem to cut the variety by suitable linear sections! 

I would like to end by reproducing the concluding remarks in my talk on Gilles at AGCT-2019:

\begin{quotation}
Gilles Lachaud has made important and lasting contributions to mathematics, especially in the study of algebraic varieties over finite fields and linear codes. His knowledge and interests were deep and wide. When he became interested in some topic, he would usually delve deeper and spend considerable time learning many aspects of it. 
As far as I have seen, he would never be in a rush to publish quickly, but would prefer to take time and be thorough.
Besides his contributions to mathematics, Gilles was an institution builder. He helped nurture an institution like  CIRM.  Also, the continuing success of the AGCT conferences owes largely to his vision and efforts. 

Other than scientific institutes and conferences, Gilles served as the President of the French Pavilion at Auroville, near Pondicherry, India. He had read or had at least browsed through significant amount of ancient and modern Sanskrit works, including the \emph{Vedas, Upanishadas}, and the scholarly treatises of \emph{Sri Aurobindo}.

Above all, Gilles was a wonderful human being, generous, warmhearted,  and kind, always willing to help others, especially students and younger colleagues. His untimely demise last year is a great loss to our subject and the community. Personally, it has been a pleasure and honor to have known him. He will certainly be missed...
\end{quotation}

\medskip
 
\section{Gilles Lachaud's work}
\medskip
\centerline{Christophe Ritzenthaler}
\medskip

Among the four contributors, I am the one who met Gilles the latest, precisely in 2006. He was then director of IML (Institut de Math\'ematiques de Luminy) in Marseille and I was starting as Ma\^itre de Conf\'erences there in the ATI (Arithm\'etique et Th\'eorie de l'Information) team. I think I owe to our common interest in genus $3$ non-hyperelliptic curves that he invited me to apply for 
this position.  

In 2005, Gilles had indeed published two articles \cite{rit-lac1,rit-lac2} on the Klein quartic (the smooth quartic $X : x^3 y+y^3z+z^3x=0$)  in {\it Moscow Mathematical Journal}\footnote{From some manuscript notes I saw, it was part of a broader program to study the geometry and arithmetic of more general family of quartics with non-trivial automorphism groups.
}. As many other people (a full book ``The eightfold way'' \cite{eightfold}  was published a few years earlier on this single curve), Gilles was fascinated by the conjunction of beautiful analytic, algebraic and arithmetic geometry and group and representation theory to address this particular genus $3$ non-hyperelliptic curve. In the first and longest paper {\it Ramanujan modular forms and the Klein quartic}, he conjugates his passion for history of mathematics and research, unravelling Ramanujan's contributions to some identities between modular forms of level $7$. The relation with modular forms boils down to the fact that $X$ is isomorphic over $\mathbb{C}$ to the modular curve $X(7)$, since they are, up to isomorphisms, the unique genus $3$ curve with automorphism group isomorphic to $\textrm{PSL}_2(\mathbb{F}_7)$. He then gives four different expressions for an $L$-series which encodes the number of points of the Klein quartic over finite fields, and derives some relations between the Weil polynomial of $X \otimes \F_p$ and the one of an elliptic curve $E$ with complex multiplication by $\Z[(1+\sqrt{-7})/2]$.  In the second article {\it The Klein quartic as a cyclic group generator}, Gilles comes back  to this last problem and gives new expressions for the Weil polynomial. He does it not only for the Klein quartic but also for some of its twists $a x^3 y+b y^3z+c z^3 x=0$ with $a,b,c \in \F_p^{\times}$. This is realized as a special case of a more general result on Faddeev curves that he proves using Davenport-Hasse's method.   These results have been expanded later by Meagher and Top in \cite{rit-mea} to all twists of the Klein quartic.

When I arrived in Marseille, it was then natural that we started a collaboration on a tantalizing question of Serre on genus $3$ curves.  In his lectures {\it Rational points on curves over finite fields} at Harvard in 1985, Serre explains that a principally polarized abelian variety $A$ of dimension $g>2$ defined over a perfect field $k$, which is geometrically a Jacobian, is not necessarily a Jacobian over $k$ (unlike in dimension $1$ or $2$). The obstruction is given by a quadratic character of $\mathrm{Gal}(\bar{k}/k)$ and is called  \emph{Serre's obstruction}. When $k \subset \mathbb{C}$ and $g=3$, he speculates in his lectures (and gives more details in a letter to Top in 2003), that this obstruction can be computed in terms of the value of a certain Siegel modular form called $\chi_{18}$. At that time, it was pure magic to us as this Siegel modular form was a product of so-called theta constants that are purely transcendental objects and which only  ``see'' the curve over $\mathbb{C}$. We could not imagine them containing arithmetic information. In order to understand Serre's intuition, we decided to work out the formulae in the particular case of Ciani quartic. They form a dimension $3$ family of quartics with automorphism group containing a group isomorphic to $(\mathbb{Z}/2\mathbb{Z})^2$. This 
structure helps decomposing their Jacobians into a product of three explicit elliptic curves up to isogeny. 
In 2000, Howe, Leprevost and Poonen in \cite{rit-leprevost} had shown how, conversely,  starting from elliptic curves $E_1,E_2$ and $E_3$ defined over a field $k$ and a ``good'' $k$-rational maximal isotropic subgroup $G \subset (E_1 \times E_2 \times E_3)[2]$, the abelian threefold $A=(E_1\times E_2\times E_3)/G$ 
 is the Jacobian over $k$ of an explicit Ciani quartic
if and only if a certain expression in terms of the coefficients of the $E_i$ is a non-zero square in $k$. They had found an algebraic expression of Serre's obstruction for $A$ in this particular case! In {\it On some questions of Serre on abelian threefolds} \cite{rit-lac3}, we worked out this expression in terms of the theta constants associated to the $E_i$, then to $A$ and eventually in terms of $\chi_{18}$ as required.  This first collaboration with Gilles went very smoothly: I had the tools to manipulate theta functions, and he knew very well the geometry of this family of quartics and useful decompositions of the symplectic group from his earlier works. This result was presented at the conference organized for Gilles' 60th birthday in Tahiti in 2007. 

 More than confirming Serre's intuition for a dimension 3 family, this work made us realize the correct normalization of the Siegel modular form we should look for. We used this knowledge to work out a complete proof of Serre's formula in our next article {\it Jacobians among abelian threefolds: a formula of Klein and a question of Serre} \cite{rit-lac4}, in collaboration with Zykin (simultaneously, Meagher got a similar result in his PhD). In order to achieve that, the best formalism was to use the theory of algebraic stacks and geometric modular forms. It sounded scary at first but Gilles was confident and calm as usual and we managed to get our way through. The cherry on top of the cake was that we could make precise a formula due to Klein linking this modular form to the discriminant of plane quartics. This formula is the first beachhead that enables connecting very precisely the world of modular forms and invariants as it was done previously in genus $1$ and $2$. We used it  recently  with Lercier to give the complete dictionary in genus $3$.
 
Three years later, I was leaving for Rennes, and we did not get the chance to collaborate more. Driven by his tireless curiosity and maybe 
triggered by a conference and a winter school at CIRM in 2014 and discussions with Kohel,  
Gilles' interest moved towards arithmetic statistics.
In his  first article on the topic \cite{rit-lac5}, {\it On the distribution of the trace in the unitary symplectic group and the distribution of Frobenius}, he considers the group $\textrm{USp}_{2g}$ which corresponds to the ``generic case'' for curves and abelian varieties in Katz-Sarnak theory. He recalls the  link between the distribution associated to the trace on $\textrm{USp}_{2g}$ and the distribution of the number of points for curves or abelian varieties over finite fields. Although this is well-known to experts, this is done with great care and with the purpose of making it explicit enough for anybody who would need to do actual 
computations. When $g=2$ and $3$, he then gives some useful expressions for the distribution law in terms of hypergeometric series and related functions and draws some inspiring graphics representing the density of curves with a certain Weil polynomial. In 2018, in {\it The distribution of the trace in the compact group of type $G_2$} \cite{rit-lac6}, he works out similar formulas for the case of a group (the only compact simple Lie group of type $G_2$) related to exponential sums and introduced by Katz. In the present volume, one can find his last work on the case of $SU(3)$.

During my time in Marseille, Gilles, as director, managed to create such a pleasant working atmosphere that the years I spent there were extraordinary. Although, as I found out later, such a position is really time-consuming, his door was always open, and he never made anybody felt indebted for the facilities he provided. 
I learned a lot from working with him and I still  meditate
about his way of addressing new problems. When I consider them as ``enemies'' that should surrender after you have attacked them from every angles, Gilles never got upset by their resistance and approached them peacefully. 
 He was also a true humanist as Fran\c cois, Michael and Sudhir mentioned, with many passions outside mathematics. In his car, stuck in the traffic on our way back from Luminy some evenings, I was always happy to listen to his insights on Japanese culture, his collaborative projects with India or lately his study  with his wife Patricia of the secrets of the painting {\it The Tempest} by Giorgione. I guess we will all miss his friendly attitude towards life. 

\bibliographystyle{amsplain}
\bibliography{rod-lac,rit-lachaud,sudhir-misha}
\end{document}